\documentclass{conm-p-l}
\usepackage{amssymb}
\usepackage[all]{xy}

\newcommand{\mHtp}{\mathbin{\stackrel{\cdot}{\otimes}}}

\newcommand{\tr}{\triangleright}
\newcommand{\tl}{\triangleleft}

\newcommand{\va}{\varphi}
\newcommand{\vt}{\vartheta}

\newcommand{\bd}{\begin{document}}
\newcommand{\ed}{\end{document}}
\newcommand{\nl}{\nu\in\Lambda}
\newcommand{\xl}{X_\nu;\nl}
\newcommand{\x}{X_\nu}
\newcommand{\la}{\langle}
\newcommand{\ra}{\rangle}
\newcommand{\ch}{\mathcal H}
\newcommand{\su}{\subseteq}
\newcommand{\B}{\Bbb}
\newcommand{\e}{\varepsilon}
\newcommand{\al}{\alpha}
\newcommand{\de}{\delta}
\newcommand{\lo}{\Longleftrightarrow}
\newcommand{\cK}{\mathcal K}
\newcommand{\bb}{\mathcal B}
\newcommand{\dd}{\mathcal D}
\newcommand{\bgd}{{\bigtriangledown}}
\newcommand{\bu}{{\bigtriangleup}}
\newcommand{\cR}{\mathcal R}
\newcommand{\di}{{\diamondsuit}}
\newcommand{\f}{\mathcal F}
\newcommand{\cb}{\mathcal{CB}}
\newcommand{\wrr}{\widetilde{\mathcal R}}
\newcommand{\br}{^\bullet{\mathcal R}}
\newcommand{\en}{E_\nu}
\newcommand{\el}{{\mathcal L}}
\newcommand{\lr}{\Longrightarrow}
\newcommand{\Long}{\Longleftarrow}
\newcommand{\q}{\quad}
\newcommand{\qq}{\qquad}
\newcommand{\cd}{\cdot}
\newcommand{\fur}{f: E\to{\B R}}
\newcommand{\furo}{f_0: E_0\to{\B R}}
\newcommand{\fuc}{f: E\to{\B C}}
\newcommand{\id}{{\bf 1}}
\newcommand{\co}{{\B C}}
\newcommand{\s}{\sigma}
\newcommand{\lm}{\lambda}
\newcommand{\ii}{\infty}
\newcommand{\mma}{\mathrel{\mathop{\otimes}\limits_{A}}}
\newcommand{\mmb}{\mathrel{\mathop{\otimes}\limits_{\bb}}}

\begin{document}
\title[Tensor products in quantum functional analysis]%
{Tensor Products in Quantum Functional Analysis:\\
the Non-Matricial Approach}
\author{A. Ya. Helemskii}
\address{Faculty of Mechanics and Mathematics\\
Moscow State University\\
119992 Moscow\\
Russia}
\email{alexander@helemskii.mccme.ru,helemskii@comail.ru}
\thanks{This research was supported by the Russian Foundation for
Basic Research (grant No. 05-01-00982).}
\begin{abstract}
As is known, there exists an alternative, "non-matricial" way to present
basic notions and results of quantum functional analysis (= operator
space theory). This approach is based on considering, instead of matrix
spaces, a single space, consisting, roughly speaking,
of vectors from the initial linear space equipped with coefficients taken
from some good operator algebra.
It seems that so far there was no systematical exposition of the theory
in the framework of the non-matricial approach. We believe, however,
that in a number of topics the non-matricial approach
gives a more elegant and transparent theory.

In this paper we introduce, using only the non-matricial
language, both quantum versions of the classical (Grothendieck) projective
tensor product of normed spaces. These versions correspond to the "matricial"
Haagerup and operator-projective tensor products. We define them in terms of
the universal property with respect to some classes of bilinear operators,
corresponding to "matricial" multiplicatively bounded and completely
bounded, and then produce their explicit constructions.
Among the relevant results, we shall
show that both tensor products are actually quotient spaces of some
"genuine" projective tensor products. Moreover, the Haagerup
tensor product is itself a "genuine" projective tensor product, however not
of just normed spaces but of some normed modules.
\end{abstract}
\maketitle

This paper deals with some questions of quantum, or quantized functional
analysis (cf. the memorable lecture of Effros~\cite{Eff}). Here
we present some basic notions and results concerning tensor
products of quantum spaces ("spaces endowed with operator space structure").
The specific feature of our exposition is that we
systematically use what can be called non-matricial or non-coordinate
approach to quantum functional analysis.

This means, to speak informally, the following thing. Usually a quantum
norm (=operator space structure) on a given linear space
is introduced by simultaneous consideration of matrices of all sizes with
entries in this space. There is, however another way. Instead of these matrix
spaces, one considers a single space, consisting, roughly speaking,
of vectors from the given space equipped with coefficients taken from some
good operator algebra. Such a replacing of scalars by operators, this time
in the capacity of coefficients to our vectors, is, of course, well in line
in the general spirit of modern "quantum", or non-commutative, mathematics.

The very fact that both approaches, matricial (coordinate) and operator
(non-coordinate), give essentially equivalent results, is known. It is
clearly indicated in the book of Pisier~\cite{Pis} who demonstrates the
virtues of the non-matricial approach in a number of important questions
(see, e.g., {\it idem}, p. 40). Besides, this was well realized by
Barry Johnson,
as one can judge from his unpublished notes. The same fact is demonstrated
in the form of theorems on the equivalence of various
categories~\cite{Web,Ng}; cf. also results on representations of bimodules
over operator algebras~\cite{mag,pop}.

But what about concrete forms acquired by the main notions of the theory
exclusively in the framework of the non-coordinate approach, without an
appeal to matrix spaces? Do we obtain some new insight on the subject?
It seems that there was no systematical exposition of quantum
functional analysis from the indicated point of view.
As to already existing monographs on
the subject, even in~\cite{Pis} the matricial approach considerably prevails,
and it is the only approach  taken in~\cite{ER,Pau,BlM}.

The choice between the two approaches is, of course, the matter of taste.
(One prefers to work with tensor products of linear operators, and another
one with Kronecker products of matrices). However, we believe
that there are some topics where the non-coordinate approach
gives more elegant and transparent theory. Especially this concerns the
notions where the matricial presentation inevitably creates the whole parade
of indices and multi-indices, and first of all that of tensor product. (The
same, as it seems to us, could be said about the duality theory, but we do
not touch this topic here).

The principal aim of this paper is to introduce, using the non-coordinate
language, both quantum versions of the classical (Grothendieck) projective
tensor product of normed spaces. These versions are usually called Haagerup and
operator-projective tensor products. We define them in terms of the
universal property with respect to some classes of bilinear operators,
corresponding to "matricial" multiplicatively bounded and completely
bounded bilinear operators, and then produce their explicit constructions.
Among the relevant results, we shall
show that both tensor products are actually quotient spaces of some
"genuine" projective tensor products. Moreover, the Haagerup
tensor product is itself a "genuine" projective tensor product, however not
of just normed spaces but of some normed modules.\footnote{The author gave
talks about these results at the conferences in Athens (June 2005) and
Bordeaux (July 2005). After the second talk Dr.~C.~-K.~Ng, who was present,
kindly informed the author that the results are essentially known to him,
and that the similar things are contained in his yet unpublished preprints.}

>From the huge number of substantial examples of quantum spaces we use here
only one pair of very illustrative twin spaces, the so-called column and
row Hilbertians. Their habits and, in particular, their behaviour as tensor
factors are well known in the matricial exposition; see, e.g.,
~\cite[Section 9.3]{ER}. We just want to show that these things,
providing the adornment of general results, well fit
in the non-coordinate presentation.

\section{Preparing the stage}

Throughout the paper, the terms operator and bioperator mean always
respectively linear operator and bilinear operator; the terms functional
and bifunctional have the similar meaning.
If $E$ and $F$ are normed spaces, then
${\mathcal B}(E,F)$  and ${\mathcal K}(E,F)$ denote, as usual,
the space of all bounded, respectively compact, operators between
these spaces, ${\mathcal B}(E)$ means ${\mathcal B}(E,E)$ and ${\mathcal K}(E)$ means
${\mathcal K}(E,E)$. The identity operator on $E$ is denoted by
${\id}_E$ or, if it is safe, just by ${\id}$.

The term {\it operator space} will be used for an arbitrary (not
necessarily closed)
subspace in ${\mathcal B}(H,K)$ for some {\it Hilbert spaces} $H$ and $K$.
So far we do not equip operator spaces with any additional structure save
induced (= operator) norm. The symbol $\mHtp$ denotes, according to the
sense, one of the following
three things. Namely,
$H\mHtp K$ is the Hilbert tensor product of respective
Hilbert spaces
whereas $a\mHtp b$ is the Hilbert tensor product of respective bounded
operators acting between
Hilbert spaces (see, e.g.,\cite[Ch.2, \S8]{baf}). Finally, if
$E\subseteq\bb(H_1,K_1)$ and
$F\subseteq\bb(H_2,K_2)$ are operator spaces, then $E\mHtp F$ is the
so-called (non-completed) {\it
spatial tensor product} of $E$ and $F$, that is the operator space
$span\{a\mHtp b;a\in E,b\in F\}\subseteq\bb(H_1\mHtp H_2,K_1\mHtp K_2)$.

If $H$ is a Hilbert space and $\xi,\eta\in H$, we denote by
$\xi\bigcirc\eta:H\to H$ the rank one operator taking
$\zeta$ to $\la\zeta,\eta\ra\xi$. Let us distinguish the obvious equalities
$$
(\xi\bigcirc\eta)(\xi'\bigcirc\eta')=
\la\xi',\eta\ra(\xi\bigcirc\eta'),
a(\xi\bigcirc\eta)=a(\xi)\bigcirc\eta\q {\rm and} \q
\|\xi\bigcirc\eta\|=\|\xi\|\|\eta\|, \eqno{(1)}
$$
where $\xi,\xi'\eta,\eta'\in H,a\in\bb(H)$.

For a short time, consider an arbitrary unital $\bb$-bimodule $X$. Take
$u\in X$. We call every projection (= self-adjoint
idempotent) $P\in{\mathcal B}$ a {\it left} (respectively, {\it right})
{\it support} of the element
$u$, if $P\cdot u=u$ (respectively, $u\cdot P=u$.
If we have both equalities, we speak about (just) a {\it support} of $u$.

Let $\|\cdot\|$ be a semi-norm on $X$.
We say that it satisfies the {\it first axiom of Ruan}
(briefly, $(RI)$) if, for every $a\in{\mathcal B}$ and $u\in X$ we have
$\|a\cd u\|,\|u\cdot a\|\le\|a\|\|u\|$.
(In the usual language of the theory of Banach algebras this means exactly
that $X$ is a contractive, or linked semi-normed $\bb$-bimodule).
We distinguish the obvious

\medskip

{\bf Proposition 1.}  {\it Suppose that the semi-norm on $X$ satisfies $(RI)$.
Then for every $u\in X$, isometric operator $a\in\bb$ and every coisometric
operator $b\in\bb$ we have $\|a\cd u\|=\|u\cd b\|=\|u\|$.}

\medskip

Further, we say that a semi-norm $\|\cdot\|$ in $X$
satisfies the {\it second axiom of Ruan} (briefly, $(RII)$) if,
whenever $u$ (respectively, $v$) in $X$ has a support $P$ (respectively, $Q$),
and these projections are orthogonal (i.e. $PQ=0$), we have
$\|u+v\|=\max\{\|u\|,\|v\|\}$. (We choose both
axioms as "non-coordinate" versions of the
known Ruan axioms for matrix-norms (cf., e.g.,~\cite[p.20]{ER} or
~\cite[p.180-181]{Pau}).

\medskip

{\bf Proposition 2}. {\it Let a semi-norm $\|\cdot\|$ on $X$ satisfy
both axioms
of Ruan, and elements $u_k\in X; k=1,...,n$ have pairwise orthogonal left
supports or pairwise orthogonal right supports. Then
$\|u_1+...+u_n\|\le(\|u_1\|^2+...+\|u_n\|^2)^{\frac{1}{2}}$.}

\medskip

$\tl$ Obviously, we can assume that $\|u_k\|\ne0$ for all $k$.
Let us concentrate on the case of left supports, say $P_k$.
At first we shall show that the assertion is true if we assume,
in addition, that our supports, as operators on $L$,
have infinite-dimensional images.

The latter condition enables us to choose isometric operators
$S_k\in\bb;k=1,...,n$ such that $S_kS_k^*=P_k$.
Consider, for every $k$, $v_k:=\frac{1}{\|u_k\|}u_k\cd S_k^*\in\cK E$. It is
easy to see that $P_k$ is a (two-sided) support of $v_k$. Moreover,
by Proposition 1, we have $\|v_k\|=1$ for all $k$. Consequently, the axiom
$(RII)$, extended by obvious way to the case of $n$ elements,
gives $\|v_1+...+v_n\|=1$. Now put $a_k:=\|u_k\|S_k\in\bb$. Then we have
$$
(v_1+...+v_n)\cd(a_1+...+a_n)=
\sum_{k,l=1}^n[\frac{1}{\|u_k\|}u_k\cd S_k^*]\cd\|u_l\|S_l=
$$
$$
\sum_{k,l=1}^n[\frac{1}{\|u_k\|}\|u_l\|u_k]\cd S_k^*S_l=u_1+...+u_n.
$$
Therefore, again by $(RI)$, we have $\|u_1+...+u_n\|\le
\|a_1+...+a_n\|\|v_1+...+v_n\|
=\|a_1+...+a_n\|$. But routine calculations, using the $C^*$-identity,
show that the norm of the
operator $a_1+...+a_n$ is exactly $(\|u_1\|^2+...+\|u_n\|^2)^{\frac{1}{2}}$.
The desired estimation follows.

In the case of arbitrary left supports $P_k$ we use the following
device. Let $nL$ be a Hilbert sum of $n$ copies $L_k$ of $L$, and
$Q_k:nL\to nL$ a projection onto $L_k$.
Obviously, there exists an isometric operator $R:L\to nL$, mapping,
for every $k$, $Im(P_k)$ into $L_k$. Choose an arbitrary isometric
isomorphism $U:nL\to L$ and note that the operator $UR\in\bb$ is isometric.

Now put $u_k':=UR\cd u_k;k=1,...,n$. Of course, operators
$P_k':=UQ_kU^*\in\bb;k=1,...,n$ are pairwise orthogonal
 projections with infinite-dimensional images. Moreover, we have
$$
P_k'\cd u_k'=UQ_kU^*\cd(UR\cd u_k)=(UQ_kU^*UR)\cd u_k=(UQ_kR)\cd u_k.
$$
But the choice of $R$ obviously implies, for all $k$, that $Q_kRP_k=RP_k$.
Therefore
$$
P_k'\cd u_k'=(UQ_kR)\cd(P_k\cd u_k)=(UQ_kRP_k)\cd u_k=UR\cd(P_k\cd u_k)=
UR\cd u_k=u_k'.
$$

Thus we find ourselves in the situation, where the
desired inequality, however with $u_k'$ instead of $u_k$, is already
established. But Proposition 1 immediately gives $\|u_k'\|=\|u_k\|$ and
$\|\sum_{k=1}^nu_k'\|=\|\sum_{k=1}^nu_k\|$. $\tr$

We proved the part of the assertion, concerning left supports. The
analogous argument provides the part, concerning right supports.
%(Alternatively, one can reduce the "right" case to the "left" one, with
%the help of a new $\bb$-bimodule structure in $X$, namely, putting
%$a\bullet x:=x\cd a^*$ and $x\bullet a:=a^*\cd x)$. $\tr$

\medskip

{\bf Remark.} This proposition, at least for normed bimodules, could be
obtained as a corollary to some deep results about representations
of some bimodules over $C^*$-algebras in terms of the so-called operator
convexity (see., e.g.,~\cite{mag,pop}). But, since our aims are
quite different, we do not need any of such a strong medicine here.

\medskip

Now we leave general $\bb$-bimodules. Choose an arbitrary
separable infinite-dimensional Hilbert space, denote it, say, by $L$ and
fix it throughout the whole scope of this paper. The operator algebras
$\bb(L)$ and $\bb(L)$ we denote, for brevity, by $\mathcal B$ and $\mathcal K$.
Instead of  $\id_L$ we shall always write just $\id$.

Let $E$ be a linear space.
Denote, again for brevity, the algebraic tensor product
${\mathcal K}\otimes E$ by ${\mathcal K}E$ and call this space the
{\it amplification} of $E$.
This is, to speak informally, "the space of formal linear combinations
of vectors from $E$ with
operator coefficients from $\mathcal K$". Accordingly, we denote an
elementary tensor $a\otimes x; a\in{\mathcal K}, x\in E$ just by $ax$.
Observe that the space $\cK E$ is a ${\mathcal B}$-bimodule
with respect to the outer multiplications, well
defined by the equalities $a\cdot bx=(ab)x$ and $bx\cdot a=(ba)x;
a\in{\mathcal B},b\in{\mathcal K},x\in E$.

In the following two propositions we consider a $\bb$-bimodule of the form
$\cK E$ and we suppose that it is equipped with a semi-norm $\|\cdot\|$
satisfying $(RI)$.

\medskip

{\bf Proposition 3}. {\it Assume that $\lim_{n\to\ii}a_n=0$ for some
$a_n\in\cK$. Then, for every
%quantum space $E$ and
$x\in E$, we have $\lim_{n\to\ii}a_nx=0$.}

\smallskip

$\tl$ As is well known, there exist $b_n,c\in\cK$ such that $a_n=b_nc;
n=1,2,...$ and
$\lim_{n\to\ii}b_n=0$ (see, e.g.,~\cite[\S11, Cor.12]{BD}).
Hence, by $(RI)$, $\|a_nx\|=\|b_n\cd(cx)\|\le
\|b_n\|\|cx\|$. The rest is clear. $\tr$

\medskip

{\bf Proposition 4}. {\it Assume that, for some projection
$p\in\cK$ of rank 1 the restriction of the given semi-norm on $\cK E$ to
the subspace $\{px; x\in E\}$ is a norm. Then the given semi-norm is
itself a norm.}

\smallskip

$\tl$ Take a non-zero element $u\in\cK E$. Then, as is well known
(see, e.g.,~\cite[Proposition 2.7.1]{baf}), it can be
represented as $\sum_{k=1}^na_kx_k$, where $a_k;k=1,...,n$ is a linearly
independent system of
compact operators, and $x_1\ne0$. Consider the system of all vector
functionals on $\cK$, that is
of those acting as $a\mapsto\la a\xi,\eta\ra;\xi,\eta\in H$. Of course,
the system of vector functionals is sufficient, i.e. for every $a\ne0$ there
exists a vector functional with a non-zero value on $a$. According to the
known property of linear span of sufficient
systems (see, e.g.,~\cite[Proposition 4.2.3]{baf}), there exist
$\xi_l,\eta_l\in H; l=1,...,m$ such that
$\sum_{l=1}^m\la a_k\xi_l,\eta_l\ra$ is 1 when $k=1$ and is 0 otherwise.

Now take a normed vector $e\in Im(p)$; obviously, $p=e\bigcirc e$. Consider
in $\cK E$ the element
$v:=\sum_{l=1}^m(e\bigcirc\eta_l)\cd u\cd(\xi_l\bigcirc e)$ and recall the
equalities (1). We see that
$$
v=\sum_{l=1}^m(e\bigcirc\eta_l)\cd(\sum_{k=1}^na_kx_k)\cd(\xi_l\bigcirc e)=
\sum_{k,l}[(e\bigcirc\eta_l)a_k(\xi_l\bigcirc e)]x_k=
$$
$$
\sum_{k,l}[\la a_k\xi_l,\eta_l\ra p]x_k=
[\sum_{l=1}^m\la a_1\xi_l,\eta_l\ra p]x_1=px_1.
$$
Therefore, by our assumption, $\|v\|\ne0$. At the
same time the axiom $(RI)$ and the triangle inequality for semi-norms give
$$
\|v\|\le\sum_{l=1}^m\|(e\bigcirc\eta_l)\cd u\cd(\xi_l\bigcirc e)\|\le
\sum_{l=1}^m\|e\bigcirc\eta_l\|\|u\|\|\xi_l\bigcirc e\|.
$$
This, of course, implies $\|u\|>0$. $\tr$

\section{Quantum spaces and completely bounded operators}

The main concepts of quantum functional analysis, in its non-coordinate
presentation, are those given in the following definition and in still
more important Definition 2.

\medskip

{\bf Definition 1}.  A {\it quantum} norm on a linear space $E$ is an
arbitrary norm on the $\bb$-bimodule ${\mathcal K}E$, satisfying both of Ruan
axioms. A {\it quantum normed space} or just a {\it quantum space} is a
linear space, equipped with a quantum norm.

(We emphasize that a quantum norm on $E$ is a (usual) norm not on $E$
itself, but on the "larger" space ${\mathcal K}E$).

\medskip

A quantum normed space becomes a "classical" normed space, if, for $x\in E$,
we put $\|x\|:=\|px\|$, where $p$ is an arbitrary projection of rank 1 on $L$.
It easily follows from Proposition 1 that
this is indeed a norm, not depending on a choice of $p$. (In fact, this norm
does not change if we replace $p$ by an arbitrary $a\in\cK;\|a\|=1$.
But we do not need this fact here). The resulting
normed space, often denoted by $\Box E$,
is called the {\it underlying normed space} of the quantum
space $E$. As to the initial quantum space, we call it
a {\it quantization} of its underlying space $\Box E$, and we call its
quantum norm a {\it quantization} of the "usual" norm on $\Box E$.
We shall see that the same, up to isometric
isomorphism, normed space can have a lot of profoundly different
quantizations. However the simplest normed space, the complex plane $\co$,
has a unique quantization.
Namely, it easily follows from axioms of Ruan that the operator norm on
$\cK\co=\cK$ is the only quantization of the norm on $\co$.

Let $F$ be a linear subspace of a quantum space $E$. Then $F$
becomes itself a quantum space with
respect to the norm on $\cK F$, well defined by $\|u\|:=\|{\id}_\cK\otimes i\|$,
where $i:F\to E$ is the natural embedding.
In this situation we say that $F$ is a {\it quantum subspace} of $E$.

We turn to a general construction providing the principal class of quantum
spaces.

Assume that a linear space $E$ is given together with an
injective operator $I:E\to{\mathcal B}(H,K)$ for some Hilbert spaces $H$ and $K$.
Clearly, $E$ becomes a normed space with respect to the induced norm, and
$I$ becomes an isometric operator which provides the identification of this
normed space with the operator space $I(E)$. But more can be said in this
situation. Consider the
operator $J:\cK E\to{\mathcal B}(L\mHtp H,L\mHtp K)$, associated with the
bioperator
$\cK\times E\to{\mathcal B}(L\mHtp H,L\mHtp K):(a,x)\mapsto a\mHtp I(x)$.
It is well known (and easy to
check) that $J$ is injective. Therefore we can endow $\cK E$ with the
respective induced norm, thus
identifying $\cK E$ with the operator space $J(\cK E)=\cK\mHtp E$.
It is easy to verify that this norm on $\cK E$ is a quantum norm on $E$,
moreover a quantization of the usual norm on the latter space.
This quantum norm on $E$, as well as respective quantum space, are called
{\it concrete} quantum norm or, respectively, quantum space
({\it associated to the injective operator $I$}, if we want to be precise).

\medskip

{\bf Remark}. As a matter of fact,
every ("abstract") quantum norm on a linear space is a concrete quantum norm
(associated to
some $I$). This is the famous Ruan Theorem (see, e.g.,~\cite[p.33]{ER}), or,
more accurately, its non-coordinate version.
However, we do not need this deep theorem here.

\medskip

If a linear space $E$ is already presented as an operator space, we always
take as $I$ the
respective natural embedding and call the resulting quantum norm and
quantum space {\it standard}. The term "standard quantization" of an
operator space or of its
norm has the similar meaning. We see that in the indicated case we just
identify $\cK E$ with the operator space $\cK\mHtp E$.

Let us distinguish two important particular cases of the concrete
quantization that will provide instructive illustrations to our future
quantum tensor products.

Take a Hilbert space, say, $H$. Consider linear (and isometric) isomorphism
$I_c:H\to{\mathcal B}(\co,H)$, taking
$x$ to the operator $1\mapsto x$, as well as another linear (and isometric)
isomorphism
$I_r:H\to{\mathcal B}({\overline H},\co)$, taking $x$ to the functional
$y\mapsto\la x,y\ra$. (Here and
thereafter ${\overline H}$ denote the complex-conjugate space of $H$).
Endow $H$ with two concrete quantum norms associated respectively with
$I_c$ and $I_r$,
and denote the resulting quantum spaces by $H_c$ and $H_r$. Obviously
the underlying normed space of both $H_c$ and $H_r$ is $H$ with its
original norm.

The quantum space $H_c$ (respectively, $H_r$) is
called the {\it column} (respectively, {\it raw}) quantization of the
Hilbert space $H$, or, if $H$
is fixed, the {\it column} (respectively, {\it raw}) Hilbertian.

The following observation considerably facilitates the work with these
Hilbertians.

\medskip

{\bf Proposition 5}. {\it Let $H$ be a Hilbert space,
$E$ is a linear space. Then every
element $w$ in $E\otimes H$ has the form $\sum_{k=1}^nx_k\otimes e_k$,
where $e_k$ is an orthonormal system in $H$, and $x_k\in E$. Moreover, if
$E$ is an operator space, and $H$ is identified with the operator space
$I_c(H)$, respectively, $I_r(H)$, then $w$, being considered in the
operator space $E\mHtp H$, has the norm}
$$
\|w\|=\|\sum_{k=1}^nx_k^*x_k\|^{\frac{1}{2}},\q{\rm respectively,}\q
\|w\|=\|\sum_{k=1}^nx_k^*x_k\|^{\frac{1}{2}}.
$$

\smallskip

$\tl$ To obtain the indicated representation, we take an arbitrary
representation, say,
$\sum_{k=1}^my_k\otimes\xi_k$, of $w$, take an orthonormal basis $e_k$ in
$span\{\xi_1,..,\xi_m\}\subset H$ and use the bilinearity of the
operation "$\otimes$".

To compute $\|w\|$ in the "column" case, we note the following.
Our $e_k$, now operators in
$\bb(\co,H)$, satisfy $e_k^*e_l=\delta^k_l{\id}_\co$. (Here and thereafter
$\delta$ is the symbol of Kronecker). Combining
this with the operator $C^*$-identity $\|w\|=\|w^*w\|^{\frac{1}{2}}$, we
easily get the desired expression.

Similar argument works in the "raw" case as well. The only modification is
that now we have
$e_ke_l^*=\delta^k_l{\id}_\co$ and
use the $C^*$-identity in the form $\|w\|=(\|ww^*\|)^{\frac{1}{2}}$. $\tr$

\medskip

Consider an important particular case of the obtained equalities.
As usual, for a given partial isometry, say $S$, on some Hilbert space, the operator $S^*S$ will
be called its initial, and $SS^*$ its final projection.

Now let $q_k\in\cK;k=1,...,n$ be arbitrary (of course, finite-dimensional)
partial isometries in $\cK$ with the same initial projection $P$
and with pairwise orthogonal final
projections. Besides, let $e_1,...,e_n$ be an orthonormal system in $H$. Put
$$
\omega:=\sum_{k=1}^nq_k^*e_k\in\cK H\q {\rm and}\q
\varpi:=\sum_{k=1}^nq_ke_k\in\cK H.\eqno{(2)}
$$

\medskip

{\bf Proposition 6}. {\t If we consider the quantum space $H_c$, then
$\|\omega\|=1$ whereas $\|\varpi\|=\sqrt{n}$. At the same time, if we
consider $H_r$, then $\|\omega\|=\sqrt{n}$ whereas $\|\varpi\|=1$.}

\smallskip

$\tl$ Take $\cK$ as $E$  and do what is prescribed by Proposition 5. Then
we see that in the "column" case $\|\omega\|^2$ is the norm of the operator
$\sum_{k=1}^nq_kq_k^*$ whereas $\|\varpi\|^2$ is that of
$\sum_{k=1}^nq_k^*q_k=nP$. The assertion in the column case immediately
follows. A similar argument establishes the "raw" case. $\tr$

\bigskip

After the introducing, by Definition 1, a certain additional structure on
linear spaces, we naturally
proceed to the discussion of maps, reacting in an appropriate way to this
structure.

Let $\va:E\to F$ be an operator between linear spaces. The operator
${\id}_\cK\otimes\va:\cK E\to\cK F$, denoted for brevity by $\va_\ii$, is called
the {\it amplification} of $\va$. Note that $\va_\ii$ is
a morphism of $\bb$-bimodules (cf. above).

\medskip

{\bf Definition 2}. Let $E$ and $F$ be quantum spaces.
The operator $\va:E\to F$ is called
{\it completely bounded} if its amplification $\va_\ii$
is a bounded operator (with respect to the relevant
quantum norms). The operator norm of $\va_\ii$ is called {\it completely
bounded norm} of $\va$ and
is denoted by $\|\va\|_{cb}$. Further, the operator $\va$ is called
{\it completely contractive}
if $\va_\ii$ is contractive (i.e. $\|\va\|_{cb}\le1$),
{\it completely isometric} if $\va_\ii$ is isometric and
{\it completely isometric isomorphism} if $\va_\ii$ is an isometric
isomorphism.

\medskip

If an operator $\va:E\to F$ between quantum spaces is bounded as an
operator between the respective
underlying normed spaces, we say that it is (just) bounded. Taking an
arbitrary rank 1
projection $p$ and passing from $\va_\ii$ to its birestriction, which
acts between $\{px; x\in E\}$ and $\{py; y\in F\}$, we
see that {\it every completely bounded operator is bounded}, and
$\|\va\|\le\|\va\|_{cb}$.

In a number of important situations the converse is also true. We need
here only one result of that kind (cf. the "matricial" Corollary 2.2.3
in~\cite{ER}).

\medskip

{\bf Proposition 7}. {\it Let $f:E\to{\B C}$ be a bounded functional on a
quantum  space. Then it is (automatically) completely bounded,
and $\|f\|_{cb}=\|f\|$}.

\smallskip

$\tl$ Consider $f_\ii:\cK E\to\cK{\B C}=\cK$ and take $u\in\cK E$. By virtue of
properties of the operator norm, we have
$$
\|f_\ii(u)\|=\sup\{|\la f_\ii(u)\xi,\eta\ra|; \xi,\eta\in L,
\|\xi\|,\|\eta\|\le1\}.
$$
Fix an arbitrary normed vector $e\in L$ and take the projection
$p=e\bigcirc e$ onto its linear span. Using the first and the second of the
equalities (1) and then the morphism property of $f_\ii$ (see above),
we have
$$
\la[f_\ii(u)](\xi),\eta\ra p=\la[f_\ii(u)](\xi),\eta\ra(e\bigcirc e)=
$$
$$
(e\bigcirc\eta)([f_\ii(u)](\xi)\bigcirc e)=
(e\bigcirc\eta)f_\ii(u)(\xi\bigcirc e)=
f_\ii[(e\bigcirc\eta)\cdot u\cdot(\xi\bigcirc e)].
$$
Therefore $|\la[f_\ii(u)](\xi),\eta\ra|=
\|f_\ii[(e\bigcirc\eta)\cdot u\cdot(\xi\bigcirc e)]\|$.

Now observe that $(e\bigcirc\eta)\cdot u\cdot(\xi\bigcirc e)$ is an
elementary tensor of the form
$px_{\xi,\eta}$ for some $x_{\xi,\eta}\in E$. (Obviously, it is the case
when $u$ is an elementary tensor, and hence it is true for all $u$).
Besides, it follows from $(RI)$ and from the third equality in (1)
that $\|x_{\xi,\eta}\|=\|x_{\xi,\eta}p\|\le
\|e\bigcirc\eta\|\|u\|\|\xi\bigcirc e\|\le\|u\|$ wherever
$\|\xi\|,\|\eta\|\le1$. Hence for the same $\xi,\eta$ we have
$$
|\la[f_\ii(u)](\xi),\eta\ra|=\|f_\ii(px_{\xi,\eta})\|=
\|f(x_{\xi,\eta})p\|=|f(x_{\xi,\eta})|\le\|f\|\|x_{\xi,\eta}\|\le\|f\|\|u\|.
$$
Taking the respective supremum, we see that $\|f_\ii\|\le\|f\|$. The rest is
clear. $\tr$

\medskip

However, the "usual" boundedness, generally speaking, does not imply  the
complete boundedness, and this is a fundamental observation of the whole
theory. Probably, the simplest and most illuminating counter-example is
provided by the identity operator ${\id}:H_c\to H_r$
where $H$ is an infinite-dimensional Hilbert space. Indeed,
by virtue of Proposition 6, for every $n$ one can find an element in $\cK H$
such that the amplification ${\id}_\ii:\cK H_r\to\cK H_c$ increases its norm
exactly in $\sqrt{n}$ times.
Thus ${\id}_\ii$ is not bounded and hence the original operator,
being "on the level of underlying normed spaces" even isometric, is not
completely bounded. The same, with obvious modifications, can be said
about the operator ${\id}:H_r\to H_c$.

\section{Completely bounded bilinear operators}

As is known, there is a universal consent in the classical functional
analysis concerning what to call bounded bioperator between normed spaces.
As to quantum functional analysis, the experience of last 15 years has
shown that there exist at least two versions of the notion of completely
bounded bioperator, each with its own advantages. We begin with the earlier
version, discovered (in the "matricial" presentation) by Christensen and
Sinclair~\cite{ChS}, 1987.

Let ${\mathcal R}:E\times F\to G$ be a bioperator, connecting three linear
spaces. Consider the bioperator ${\mathcal R}_s:\cK E\times\cK F\to\cK G$,
associated with the 4-linear operator $\cK\times E\times\cK\times F\to\cK G:
(a,x,b,y)\mapsto ab\cR(x,y)$. (Otherwise, ${\mathcal R}_s$ is
well-defined by taking a pair $(ax,by)$ to $ab\cR(x,y)$). This bioperator is
called the {\it strong amplification} of $\cR$.

\medskip

{\bf Definition 3}. Let $E,F$ and $G$ be quantum spaces. A bioperator
${\mathcal R}:E\times F\to G$
is called {\it strongly completely bounded}\footnote{In the pioneering
paper~\cite{ChS} and in a lot of other papers and books, up to the present
time, such a bioperator (or, more precisely, its matricial version)
is called just completely bounded. However, in some
other books and papers, notably in the influential textbook of
Effros and Ruan~\cite{ER}, it is called multiplicatively bounded whereas the
term "completely bounded" is used for the "matricial prototype" of what we
call here weakly completely bounded bioperator.}if its strong amplification
${\mathcal R}_s$ is a bounded bioperator (with respect to the relevant
quantum norms). The bioperator norm of ${\mathcal R}_s$ is called {\it strong
completely bounded norm} of $\cR$ and is denoted by $\|\cR\|_{scb}$. Further,
the bioperator $\cR$ is called {\it strongly completely contractive} if
${\mathcal R}_s$ is contractive (i.e. $\|\cR\|_{scb}\le1$).

\medskip

In order  to introduce another version of complete boundedness for
bioperators, we need some preparation. We would like to have an operation
that imitates the tensor multiplication of operators on our fixed
Hilbert space $L$, but does not lead out of this space.

By virtue of Fischer-Riesz Theorem, there exists a unitary isomorphism
$\iota: L\to L\mHtp L$. Take one and fix it throughout this paper. (It does
not matter which one we choose). Our $\iota$ gives rise to the isometric
$^*$-isomorphism $\varkappa:={\mathcal B}(L\mHtp L)\to{\mathcal B}:
a\mapsto\iota^* a\iota$.

Let us use, for the operator $\varkappa(a\mHtp b)\in{\mathcal B};
a,b\in{\mathcal B}$, the brief notation
$a\di b$. Obviously, we have the identities
$$
(a\di b)(c\di d)=ac\di bd,(a\di b)^*=a^*\di b^*\q {\rm and} \q
\|a\di b\|=\|a\|\|b\|.\eqno{(3)}
$$
Besides,  $a,b\in\cK$ implies $a\di b\in\cK$.

Now let ${\mathcal R}$ be as above. Consider the bioperator
${\mathcal R}_w:\cK E\times\cK F\to\cK G$, associated with the 4-linear operator
$\cK\times E\times\cK\times F\to\cK G:(a,x,b,y)\mapsto (a\di b)\cR(x,y)$
(and well-defined by taking $(ax,by)$ to $(a\di b)\cR(x,y)$). This
bioperator is called the {\it weak amplification} of $\cR$.

\medskip

{\bf Definition 4}. Let $E,F$ and $G$ be quantum spaces. A bioperator
${\mathcal R}:E\times F\to G$ is called {\it weakly completely bounded}
if its weak amplification is a bounded bioperator. The bioperator norm of
${\mathcal R}_w$ is called
{\it weak completely bounded norm} of $\cR$ and is denoted by $\|\cR\|_{wcb}$.
The bioperator $\cR$
is called {\it weakly completely contractive} if ${\mathcal R}_w$ is contractive.

\medskip

Now let us widen the field of action of the operation "diamond".
Namely, for a linear space $E$ and $a\in\cK$
we consider the operators $_a\di,\di_a:\cK E\to\cK E$, associated
with the bioperators $\cK\times E\to\cK E$ taking $(b,x)$ respectively to
$(a\di b)x$ and $(b\di a)x$. Then for $a\in\cK$ and $u\in\cK E$ we put
$a\di u:=_a\di(u)$ and $u\di a:=\di_a(u)$. Obviously, both new "diamond
multiplications" are uniquely determined by their bilinearity and the
equations
$$
a\di bx=(a\di b)x,\quad {\rm respectively}\quad bx\di a=(b\di a)x;\quad
a,b\in\cK,x\in E.
$$
Mention the useful formulae
$$
(a\di b)\cd(c\di u)=ac\di(b\cd u),\quad (a\di u)\cd(b\di c)=ab\di(u\cd c),
$$
$$
(a\di b)\cd(u\di c)=(a\cd u)\di bc\q{\rm and}\q
(u\di a)\cd(b\di c)=(u\cd b)\di ac,\eqno{(4)}
$$
where $u\in\cK E$, and other letters denote operators in $\cK$ or, if it
is sensible, in $\bb$. (With the help of the first equality in (3), one can
easily check them for elementary tensors and then use the bilinearity).

\medskip

{\bf Proposition 8}. {\it Let $E$ be a quantum space and $P\in\cK$ a
projection of finite rank. Then,
for every $u\in\cK E$, we have $\|P\Diamond u\|=\|u\Diamond P\|=\|u\|$.}

\smallskip

$\tl$ Let us begin with a projection of rank one, say, $p$. Fix a normed
vector, say $e$, in its image and consider the isometric operator
$\rho:L\to L\mHtp L:\xi\mapsto e\otimes\xi$. Since
$\rho^*$ is uniquely determined by the taking $e\otimes\xi$ to $\xi$ and
$e'\otimes\xi$ to 0 for all $e':e'\perp e$, we easily see that
$\rho a\rho^*=p\mHtp a$ for all $a\in\bb$. Therefore
if we introduce the isometric operator $S_p:=\iota^*\rho\in\bb$, we have
$$
S_paS_p^*=\iota^*\rho a \rho^*\iota=\iota^*(p\mHtp a)\iota=
\varkappa(p\mHtp a)=p\di a.
$$
Consequently, we have $p\di u=S_p\cd u\cd S_p^*$ for all elementary tensors
in $\cK E$ and hence, by bilinearity, for all $u\in\cK E$. Proposition 1
immediately  implies $\|p\di u\|=\|u\|$.

Now let $P$ be a projection of rank $N$ on $L$. Then, for some pairwise
orthogonal
projections $p_1,...,p_N$ of rank one, we have $P=\sum_{k=1}^Np_k$. Take
$u\in\cK E$. Then $P\di u=\sum_{k=1}^Np_k\di u$, and elements
$p_k\di u$ have pairwise orthogonal supports, namely $p_k\di{\id}$. Therefore
the obvious extension of $(RII)$ to the case of
several elements gives $\|P\di u\|=\max\{\|p_k\di u\|;k=1,...,n\}$. Hence
$\|P\di u\|=\|u\|$. A similar argument provides $\|u\di P\|=\|u\|$. $\tr$

\medskip

{\bf Theorem 1.} (cf. the "matricial prototype" in~\cite[p.150]{ER}) {\it
Let ${\mathcal R}:E\times F\to G$ be a strongly completely bounded bioperator
between quantum spaces. Then ${\mathcal R}$ is weakly completely bounded, and}
$\|{\mathcal R}\|_{wcb}\le\|{\mathcal R}\|_{scb}$.

\smallskip

$\tl$ Take a projection $P$ of finite rank on $L$. For elementary tensors
$ax\in\cK E$ and $by\in\cK F$ we have
$$
\cR_w([ax]\cd P,P\cd[by])=\cR_w([aP]x,[Pb]y)=(aP\di Pb)\cR(x,y)=
$$
$$
(a\di P)(P\di b)\cR(x,y)=\cR_s([a\di P]x,[P\di b]y)=\cR_s(ax\di P,P\di by).
$$
Therefore, by bilinearity, we have $\cR_w(u\cd P,P\cd v)=\cR_s(u\di P,P\di v)$
for all $u\in\cK E$ and $v\in\cK F$. Hence
$\|\cR_w(u\cd P,P\cd v)\|\le\|\cR_s\|\|u\di P\|\|P\di v\|$ and, taking into
account the previous proposition, we have
$$
\|\cR_w(u\cd P,P\cd v)\|\le\|\cR_s\|\|u\|\|v\|.\eqno{(5)}
$$
Now take a sequence $P_n$ of finite-dimensional projections in $L$,
providing an approximate identity in $\cK$.
Obviously, we have $\lim_{n\to\ii}aP_n\di P_nb=a\di b$ for
all $a,b\in\cK$. Therefore Proposition 3 easily implies that
$\lim_{n\to\ii}\cR_w(u\cd P_n,P_n\cd v)=
\cR_w(u,v)$ in $\cK G$ for all elementary tensors $u\in\cK E, v\in\cK F$.
Hence, by bilinearity, the same is true for all
$u\in\cK E,v\in\cK F$. Combining this with (5), we have
$\|\cR_w(u,v)\|\le\|\cR_s\|\|u\|\|v\|$. The rest is clear. $\tr$

\smallskip

Similarly to what we have seen in the case of operators, a weakly (and hence
strongly) completely bounded operator $\cR$ between concrete quantum spaces
is automatically bounded as bioperator between the respective
underlying spaces, and we have $\|\cR\|\le\|\cR\|_{wcb}$.
Indeed, if $p\in\cK$ is a projection of rank one, the same is true for
$p\di p$. Therefore for $x\in E$ and $y\in F$ we have
$$
\|\cR(x,y)\|=\|(p\di p)\cR(x,y)\|=\|\cR_w(px,py)\|\le
\|\cR\|_{wcb}\|px\|\|py\|=\|\cR\|_{wcb}\|x\|\|y\|,
$$
and the desired fact follows.

Again, like in the case of operators, in a number of concrete
situations the converse is true.

\medskip

{\bf Proposition 9}. {\it Suppose that $E$ and $F$ be quantum spaces,
$f$ and $g$ are bounded functionals respectively on $E$ and $F$, and
$f\times g:E\times F\to{\B C}$ is the bifunctional, acting
as $(x,y)\mapsto f(x)g(y)$. Then $f\times g$ is strongly and hence weakly
completely bounded, and $\|f\times g\|_{scb}=\|f\times g\|_{wcb}=\|f\|\|g\|$.}

\smallskip

$\tl$  Obviously we have $\|f\times g\|=\|f\|\|g\|$ and hence
$\|f\|\|g\|\le\|(f\times g)_w\|$. Therefore, by virtue of Theorem 1, it
is sufficient to show that $\|(f\times g)_s\|\le\|f\|\|g\|$.

Taking elementary tensors and using the bilinearity, we easily see that
$(f\times g)_s:\cK E\times\cK F\to\cK$ acts as $(u,v)\mapsto f_\ii(u)g_\ii(v)$.
>From this, with the help of Proposition 7, we have
$$
\|(f\times g)_s(u,v)\|\le\|f_\ii(u)\|\|g_\ii(v)\|\le\|f\|\|g\|\|u\|\|v\|.
$$
The rest is clear. $\tr$

\medskip

The following property of weakly completely bounded bioperators, as we
shall soon see, has no "strong" analogue. For a bioperator
${\mathcal R}:E\times F\to G$, acting between linear spaces, put
${\mathcal R}^{op}:F\times E\to G:(y,x)\mapsto{\mathcal R}(x,y)$.

\medskip

{\bf Proposition 10}. {\it Suppose that ${\mathcal R}$ acts between quantum
spaces, and it is weakly completely bounded. Then ${\mathcal R}^{op}$ is also
weakly completely bounded, and $\|{\mathcal R}^{op}\|_{wcb}=\|{\mathcal R}\|_{wcb}$.}

\smallskip

$\tl$ Consider the flip operator $\bgd:L\mHtp L\to L\mHtp L$, well-defined by
$\xi\otimes\eta\mapsto\eta\otimes\xi; \xi,\eta\in L$. It gives rise to
another unitary operator, namely
$\bu:=\iota^*\bgd\iota:L\to L$. Since $\bgd(a\mHtp b)\bgd=b\mHtp a; a,b\in\bb$,
we have, for the same $a,b$ the equality $\bu(a\di b)\bu=b\di a$.

Now consider $({\mathcal R}^{op})_w:\cK F\times\cK E\to\cK G$. It easily
follows from the latter equality that we have $({\mathcal R}^{op})_w(v,u)=
\bu\cd{\mathcal R}_w(u,v)\cd\bu$ for elementary tensors
and hence, by bilinearity, for all elements $u\in\cK E,v\in\cK F$. Since
$\bu$ is a unitary, Proposition 1 gives $(\|{\mathcal R}^{op})_w(v,u)\|=
\|{\mathcal R}_w(u,v)\|$ for all $v\in\cK F,u\in\cK E$. The rest is
clear. $\tr$

\medskip

Again, column and raw Hilbertians provide several excellent illustrations.

\medskip

{\bf Proposition 11}. Every bounded bifunctional
$f:H_r\times K_c\to\co$, where $H$ and $K$ are Hilbert spaces, is
(automatically) strongly and hence weakly completely bounded. Moreover,
$$
\|{f}\|_{scb}=\|{f}\|_{wcb}=\|{f}\|.
$$

$\tl$ By virtue of the definition of the column and row Hilbertians,
elements of $\cK K_c$  are identified with operators from
$L=L\mHtp{\B C}$ into $L\mHtp K$, and, in particular, the elementary tensor
$ax$ transforms to the operator $\xi\mapsto a(\xi)\mHtp x$. At the same
time elements of $\cK H_r$ are identified with operators from
$L\mHtp{\overline H}$ into
$L=L\mHtp{\B C}$, and the elementary tensor $by$ transforms to the
operator well defined by $\eta\mHtp z\mapsto b(\eta)\la y,z\ra$. (We
emphasize that
$\la\cdot,\cdot\ra$  denotes in our argument the inner product in
$H$, and not in ${\overline H}$).

As is well known, our $f:H\times K\to\co$ gives rise to a bounded operator
$\va:{K}\to{\overline H}$, well defined by
$\la y,\va(z)\ra={f}(y,z); y\in{\overline H}, z\in K$, and we have
$\|{f}\|=\|\va\|$.
Consider, for $u\in\cK K$ and $v\in\cK{H}$ in their capacity of operators,
the diagram
$$
L\stackrel{u}{\longrightarrow}L\mHtp K
\stackrel{{\id}\mHtp\va}{\longrightarrow}L\mHtp{\overline H}
\stackrel{v}{\longrightarrow}L.
$$
If $u=ax$ and $v=by$, then the easy calculation shows that the
respective operator composition takes $\xi\in L$ to $f(y,x)ba(\xi)$, that is
to $[f_s(v,u)](\xi)$. (Here $\cK\otimes\co$, the range of $f_s$, is, of course,
identified with $\cK$). By bilinearity, we have that our composition is
$f_s(v,u)$ for all $v\in\cK H$ and $u\in\cK K$. But then
$\|{f}_s(v,u)\|\le\|v\|\|{\id}\mHtp\va\|\|u\|=\|{f}\|\|v\|\|u\|$.
The rest is clear. $\tr$

\medskip

Combining this proposition with Proposition 10, we get

\medskip

{\bf Corollary 1}. Every bounded bifunctional
$f:H_c\times K_r\to\co$, where $H$ and $K$ are Hilbert spaces, is
(automatically) weakly completely bounded, and $\|{f}\|_{wcb}=\|{f}\|$.

\medskip

But why only weakly? Now the time of counter-examples arrived. Probably,
the most transparent of them are based on the bifunctional of inner product
$\la\cdot,\cdot\ra:H\times{\overline H}\to\co$.

It is easy to see that the strong amplification of this bifunctional
takes the pair
$(\omega,\varpi)$, introduced in Section 2, to the operator
$nP$ which has, of course, the norm $n$. Since we can take an arbitrary $n$,
Proposition 6 implies that our bifunctional,
being considered on $H_c\times{\overline H}_r$ is not strongly
completely bounded. But we already know that such a bifunctional is
weakly completely bounded.
This shows, first, that the words "strong" and "weak" used here are not
for nothing, and, second,
that a bioperator $\cR$ can well be strongly completely bounded
whereas $\cR^{op}$ is not.

Finally, to display a bounded bioperator which is not even weakly
completely bounded, one can take the same bifunctional but considered on
$H_c\times{\overline H}_c$ or on $H_r\times{\overline H}_r$. Indeed,
the weak amplification of our bifunctional obviously takes the pair
$(\omega,\omega)$ to $\sum_{k=1}^nq_k^*\di q_k^*$ and
$(\varpi,\varpi)$ to $\sum_{k=1}^nq_k\di q_k$. In both cases, as is
easy to see, we get an operator of norm $\sqrt{n}$. Again, Proposition 6
immediately gives what we want.

One more example, this time of more general nature, deserves our special
attention. Suppose that
$E$ and $F$ are explicitly presented as operator spaces.
Consider their spatial tensor product $E\mHtp F$ (see Section 1)
and equip it with the standard quantum norm. This means, as we remember,
that $\cK(E\mHtp F)$ is identified with the operator space
$\cK\mHtp(E\mHtp F)$ (as well as $\cK E=\cK\mHtp E$ and $\cK F=\cK\mHtp F$).

\medskip

{\bf Proposition 12}. {\it The bioperator ${\mathcal T}:E\times F\to E\mHtp F$,
acting as $(x,y)\mapsto x\mHtp y$, is strongly completely contractive}.

\smallskip

$\tl$ Consider ${\mathcal T}_s:\cK E\times\cK F\to\cK\mHtp(E\mHtp F)$, in the
current situation acting
between $(\cK\mHtp E)\times(\cK\mHtp F)$ and $\cK\mHtp(E\mHtp F)$.
Let $E$ is presented as a subspace of $\subseteq\bb(H_1,K_1)$, and $F$ as that
of $\subseteq\bb(H_2,K_2)$. Then we have the inclusions
$E\mHtp F\subseteq\bb(H_1\mHtp H_2,K_1\mHtp K_2),
\cK\mHtp E\subseteq\bb(L\mHtp H_1,L\mHtp K_1),
\cK\mHtp F\subseteq\bb(L\mHtp H_2,L\mHtp K_2)$ and
$\cK\mHtp(E\mHtp F)\subseteq\bb(L\mHtp(H_1\mHtp H_2),L\mHtp(K_1\mHtp K_2))$.

Take $u\in\cK\mHtp E$ and $v\in\cK\mHtp F$. Introduce the operators
$U:=u\mHtp{\id}_{K_2}:L\mHtp H_1\mHtp K_2\to L\mHtp K_1\mHtp K_2$ and
$V:L\mHtp H_1\mHtp H_2\to L\mHtp H_1\mHtp K_2$ that coincides with
$v\mHtp{\id}_{H_1}$ after the natural identification of
Hilbert tensor products with factors presented in different order.
We see that the
composition $UV$ maps $L\mHtp H_1\mHtp H_2$ into  $L\mHtp K_1\mHtp K_2$.

Now assume, for a moment, that $u$ and $v$ are elementary tensors, say
$ax$ and $by$. Then
$UV$ obviously takes the elementary tensor
$\xi\otimes\eta\otimes\zeta\in L\mHtp H_1\mHtp H_2$ to
$ab(\xi)\otimes x(\eta)\otimes y(\zeta)\in L\mHtp K_1\mHtp K_2$. This
means that in the considered case $UV$, after the
installing of the proper brackets in the respective Hilbert tensor products,
is not other thing than ${\mathcal T}_s(u,v)$. It follows, by the bilinearity
of the relevant operations, that
the same is true for all $u$ and $v$. Consequently we have
$$
\|{\mathcal T}_s(u,v)\|=\|UV\|\le\|U\|\|V\|=
\|u\mHtp{\id}_{K_2}\|\|v\mHtp{\id}_{H_1}\|=\|u\|\|v\|.
$$
The rest is clear. $\tr$

\section{The Haagerup tensor product}

The role of tensor products in quantum functional analysis is even more
important than in classical functional analysis. Their raison d'\^etre is
essentially the same: they "linearize" bilinear operators. As a "classical"
prototype, both of our quantum tensor products have the projective
tensor product of Grothendieck of normed spaces (cf.,
e.g.,~\cite[Ch.2\S7]{baf}). In fact we shall show that
their constructions are slight complifications of the construction of the
projective tensor product (see our introduction).

Up to the rest of our paper, we fix arbitrary quantum spaces $E$ and $F$.

\medskip

{\bf Definition 9} (cf. "classical" Definition 2.8.3 {\it idem}).
We say that the pair $(\Theta,\theta)$, consisting of a quantum space
$\Theta$ and a strongly (respectively, weakly) completely contractive
operator $\theta:E\times F\to\Theta$, is the {\it Haagerup tensor product}
(respectively, the {\it four-named tensor product})\footnote{The origin of
the first term is explained, e.g., in~\cite[p.173]{ER}. Our "four names"
are those of Effros/Ruan and Blecher/Paulsen who have discovered the notion
(in the "matricial" presentation) simultaneously and independently
in~\cite{EfR} and~\cite{BPa},1991. This second version of the tensor
product is
called just projective tensor product in many papers and textbooks, notably
in~\cite{ER}. We feel, however, that in our "non-coordinate" presentation
these words could create some confusion with the classical meaning of the
term, due to Grothendieck. Moreover, if we compare the respective quantum
norms, the first (Haagerup) tensor product is, to speak informally,
even "more projective" than the second one (see Theorem 3 below).},
if, for every strongly (respectively, weakly) completely contractive
bioperator $\cR:E\times F\to G$, where $G$ is some third quantum space,
there exists a completely contractive operator $R:\Theta\to G$ such that
the diagram

\medskip

$$
\xymatrix{
E\times F \ar[d]^\theta \ar [dr]^\cR\\
\Theta \ar [r]^R & G      }
$$

\medskip

\noindent is commutative.

This (so far hypothetical) operator $R$ is called
{\it associated with the bioperator $\cR$}.

\medskip

{\bf Remark.} We emphasize that normed spaces considered in this paper,
are not, generally speaking, assumed to be complete. We only note that
there exists a substantial notion of a quantum {\it Banach} space, and
both mentioned types of quantum tensor products have their Banach
versions. The latter can be constructed with the help of the proper quantum
version of the classical construction of the completion of a normed space.
But we do not touch this circle of questions in this paper.

\medskip

Using the standard general-categorical argument, based on the uniqueness of
the initial object in a category, one can easily prove that, for each of
our versions of the quantum tensor product,
the relevant uniqueness theorem is valid.
Namely, if $(\Theta_k, \theta_k); k=1,2$ are two Haagerup (respectively,
four-named) tensor products of $E$ and $F$, then there exists a
completely isometrical isomorphism
$I:\Theta_1\to\Theta_2$, making the diagram

\medskip

$$
\xymatrix{
& E\times F \ar[dl]^{\theta_1} \ar [dr]^{\theta_2}\\
\Theta_1 \ar [rr]^I & & \Theta_2      }
$$

\medskip

\noindent commutative. But we shall not do it here.

We proceed to the existence theorem for both types of quantum tensor
products. This will be proved by displaying their explicit constructions.
In both cases (just as in the classical case) our $\Theta$,
as a linear space, is the algebraic tensor product $E\otimes F$, and
$\theta$, as a bioperator, is the canonical bioperator
$\vartheta:E\times F\to E\otimes F:(x,y)\mapsto x\otimes y$. Accordingly,
our task is to supply the amplification $\cK(E\otimes F)$ by two
appropriate norms. (We recall again that in this paper we consider the
"non-completed" versions of our tensor products; cf. the previous remark).

We begin with the Haagerup tensor product. Let us consider the strong
amplification $\vt_s:\cK E\times\cK F\to\cK(E\otimes F)$ of $\vt$ and use the
notation $u\odot v$ instead of $\vt_s(u,v)$. We see that the operation
$\odot$, the so-called "Effros symbol", is bilinear. (It is
the non-coordinate analogue of the known "Effros symbol" in the matricial
exposition; cf., e.g.,~\cite[\S9.1]{ER}). Thus it is
well-defined on elementary tensors by
$ax\odot by=(ab)(x\otimes y)$. Note also that, for every
$a\in\bb,u\in\cK E,v\in\cK F$ we have $u\cd a\odot v=u\odot a\cd v$. (In
other words, the bioperator $\odot$ is balanced with respect to the right
outer multiplication in $\cK E$ and the left outer multiplication in $\cK F$).

Let $\odot:\cK E\otimes\cK F\to\cK(E\otimes F)$ be the linear operator,
associated with $\vt_s$; it is well-defined by
$\odot(u\otimes v):=u\odot v$. Since every element of $\cK$ is a product
of other elements, every elementary tensor and hence arbitrary element
in $\cK(E\otimes F)$ belongs to the image of $\odot$; in other words,
$\odot$ is surjective (Soon, in Proposition 16, we shall see that the
same is true even for the map $\vt_s$ itself). Therefore $\cK(E\otimes F)$
can be identified with a quotient space of $\cK E\otimes\cK F$. Consequently,
every norm on $\cK E\otimes\cK F$ gives rise to its quotient semi-norm on
$\cK(E\otimes F)$, and the latter semi-norm is uniquely determined by the
claim that $\odot$ is a coisometric operator (i.e. it maps the open unit
ball of $\cK E\otimes\cK F$ onto that of $\cK(E\otimes F$).

Now we apply this construction to the projective norm $\|\cd\|_p$ in
$\cK E\otimes\cK F$, that is to the norm of the (non-completed) Grothendieck
projective tensor product of normed spaces $\cK E$ and $\cK F$. (We recall,
that, for given normed spaces $X$ and $Y$, and $u\in X\otimes Y$, $\|u\|_p:=
\inf\{\sum_{k=1}^n\|x_k\|\|y_k\|\}$, where the infimum is taken for all
possible representations of $u$ as $\sum_{k=1}^nx_k\otimes y_k;
x_k\in X,y_k\in Y$; cf., e.g.,~\cite[Ch.2\S7]{baf}. {\it The resulting
quotient semi-norm on $\cK(E\otimes F)$ is denoted by
$\|\cd\|_h$}.\footnote{For all evidence, the subindex "h" is used
in the literature to honour Uffe Haagerup. Well, "H" is everywhere reserved
for Hilbert...} Thus, for every $U\in\cK(E\otimes F)$ we have
$$
\|U\|_h:=\inf\{\sum_{k=1}^n\|u_k\|\|v_k\|\},\eqno{(6)}
$$
where the infimum is taken for all possible representations
of $U$ as $\sum_{k=1}^nu_k\odot v_k;u_k\in\cK E,v_k\in\cK F$.

Note that $\cK E\otimes\cK F$, being a (non-module) tensor product of a left
and  a right  $\bb$-bimodules, has the canonical
structure of a $\bb$-bimodule, and $\odot$ is obviously a morphism of
respective $\bb$-bimodules.

\medskip

{\bf Proposition 13}. {\it The semi-norm
$\|\cd\|_h$ in the $\bb$-bimodule $\cK(E\otimes F)$ satisfies the first
axiom of Ruan}.

\medskip

$\tl$ It is well known and easy to check that, for a normed algebra $A$,
the (non-module) projective tensor product of a left contractive normed
$A$-module and a right contractive normed $A$-module is a
contractive $A$-bimodule.
This concerns, in particular, the $\bb$-bimodule $(\cK E\otimes\cK F,
\|\cd\|_p)$. But it was observed that $\cK(E\otimes F)$ is the image of the
latter bimodule with respect to a coisometric
$\bb$-bimodule morphism. The rest is clear. $\tr$

\medskip

{\bf Proposition 14}. {\it Let $G$ be a quantum space, $\cR:E\times F\to G$
a strongly completely bounded bioperator, $R:E\otimes F\to G$ the
associated linear operator. Then the amplification
$R_\ii:\cK(E\otimes F)\to\cK G$ is a bounded operator with respect to the
semi-norm $\|\cd\|_h$ and the given quantum norm on $G$.
Moreover, $\|R_\ii\|=\|\cR\|_{scb}$}.

\smallskip

$\tl$ Consider the diagram

\medskip

$$
\xymatrix{
\cK E\otimes\cK F \ar[d]^\odot \ar [dr]^{R^s}\\
\cK(E\otimes F) \ar [r]^{R_\ii} & \cK G      }
$$

\medskip

\noindent where $R^s$ is the operator, associated with the strong
amplification $\cR_s:\cK E\times\cK F\to\cK G$ of the bioperator $\cR$.
By the universal property of the projective norm, we have
$\|\cR\|_{scb}=\|R^s\|$. Further, routine calculations with elementary tensors in $\cK E\otimes\cK F$
show that this diagram is commutative. It follows, taking into account that
$\odot$ is a coisometric operator, that $\|R_\ii\|=\|R^s\|$.
The rest is clear. $\tr$

\medskip

{\bf Proposition 15}. (As a matter of fact) $\|\cd\|_h$ {\it is a norm}.

\smallskip

$\tl$ Combining Propositions 4 and 13, we see that it is sufficient to
show that, for a non-zero elementary tensor $aw; a\in\cK,w\in E\otimes F$
we have $\|aw\|_h>0$. Since $w\ne0$, it is well known (cf.the proof
of~\cite[Proposition 2.7.6]{baf}) that there exist
bounded functionals $f:E\to\co$ and
$g:F\to\co$ such that $(f\otimes g)w\ne0$. Now
put in the previous proposition $\cR:=f\times g:E\times F\to\co$. By virtue
of Proposition 9,  $\cR$ is strongly completely bounded, and
$\|\cR\|_{scb}=\|f\|\|g\|$. Since in our case  $R=f\otimes g$, Proposition 14
gives $\|R_\ii\|=\|f\|\|g\|$. But, obviously, $R_\ii(aw)=
[(f\otimes g)(w)]a$. From this, since $a\ne0$ and $(f\otimes g)w\ne0$, we
have $\|R_\ii\|\|aw\|_h\ge\|R_\ii(aw)\|>0$. The rest is clear. $\tr$

\medskip

>From now on, we shall call $\|\cdot\|_h$ {\it Haagerup norm}.

Now we shall show that the expression (6) for the Haagerup norm can be
simplified. Apart from the independent interest, this will help to shorten
some further proofs.

\medskip

{\bf Proposition 16}. {\it Every $U\in\cK(E\otimes F)$ can be represented as}
(a single "Effros symbol") $u\odot v;u\in\cK E,v\in\cK F$. {\it Moreover,
we have
$$
\|U\|_h:=\inf\{\|u\|\|v\|\},\eqno{(7)}
$$
where the infimum is taken for all possible representations of $U$ in the
indicated form.}

\smallskip

$\tl$ Take $\e>0$. Because of the equality (6), there exists a
representation $U=\sum_{k=1}^nu_k\odot v_k$ such that
$\sum_{k=1}^n\|u_k\|\|v_k\|<\|U\|_h+\e$.
Choose isometric operators $S_1,...,S_n\in\bb$ with pairwise orthogonal
images
and put $u:=\sum_{k=1}^n u_k\cd S_k^*, v:=\sum_{k=1}^nS_k\cd v_k$. Then,
since the bioperator $\vt_s$  is balanced and $S_k^*S_l=\de^k_l{\id}$,
we have
$u\odot v=\sum_{k,l=1}^nu_kS_k^*S_l\odot v_l=\sum_{k=1}^nu_k\odot v_k=U$.

Further, pairwise orthogonal projections $S_kS_k^*$ are right supports of
$u_k\cd S_k^*$ and left supports of $S_k\cd v_k;k=1,...,n$. Therefore,
by virtue of Proposition 2, we have $\|u\|\le
(\sum_{k=1}^n\|u_k\|^2)^{\frac{1}{2}}$
and $\|v\|\le(\sum_{k=1}^n\|v_k\|^2)^{\frac{1}{2}}$. But using, if
necessary, scalar
multiples, we have a right to assume that $\|u_k\|=\|v_k\|$ for all $k$.
It follows that $\|u\|\|v\|\le\sum_{k=1}^n\|u_k\|^2<\|U\|_h+\e$. Since
$\e>0$ is arbitrary, the infimum
in (7) is not bigger than $\|U\|_h$. The inverse inequality is obvious. $\tr$

\medskip

{\bf Proposition 17}. {\it The Haagerup norm on $\cK(E\otimes F)$ satisfies the
second axiom of Ruan}.

\smallskip

$\tl$ Let $U,V\in\cK(E\otimes F)$ have orthogonal supports $P_1$ and $P_2$.
Obviously, we have a right to assume that $\|U\|_h>\|V\|_h$.

Using the previous proposition, take $\e$ with $0<\e<\|U\|_h-\|V\|_h$ and
representations $U=u_1\odot v_1$, $V=u_2\odot v_2$ such that
$\|u_1\|\|v_1\|<\|U\|_h+\e$ and $\|u_2\|\|v_2\|<\|V\|_h+\e$. We can, of
course, assume that $u_1=P_1\cd u_1, u_2=P_2\cd u_2,
v_1=v_1\cd P_1, v_2=v_2\cd P_2$ and also  $\|u_1\|\ge\|u_2\|$ and
$\|v_1\|\ge\|v_2\|$.
Finally, using the same device as at the end of the proof of Proposition 2,
we can assume that the images of $P_1$ and $P_2$ are infinite-dimensional.

Now take isometric operators $S_k$ with $S_kS_k^*=P_k; k=1,2$ and put
$u:=u_1\cd S_1^*+u_2\cd S_2^*$, $v:=S_1\cd v_1+S_2\cd v_2$. Routine
calculations show that $u\odot v=U+V$ and hence $\|U+V\|_h\le\|u\|\|v\|$.
But we have
$u_k\cd S_k^*=P_k\cd(u_k\cd S_k^*)\cd P_k; k=1,2$. Therefore,
because the norm on $\cK E$ satisfies $(RII)$, $\|u\|=
\max\{\|u_1\cd S_1^*\|,\|u_2\cd S_2^*\|\}$ and hence, by Proposition 1,
$\|u\|=\max\{\|u_1\|,\|u_2\|\}=\|u_1\|$. Similarly we have
$\|v\|=\|v_1\|$. Consequently $\|U+V\|_h\le\|u_1\|\|v_1\|<\|U\|_h+\e$.
Since $\e>0$ is arbitrary, we have $\|U+V\|_h\le\|U\|_h$.
Because of $P_1\cd(U+V)=U$, the inverse
inequality follows from $(RI)$. The rest is clear. $\tr$

\medskip

Combining Propositions 13 and 17, we obtain

\medskip

{\bf Corollary 2}. {\it The Haagerup norm is a quantum norm on $E\otimes F$.}

\medskip

We denote the constructed quantum space by $E\otimes_hF$. The same symbol
will denote the underlying normed space; this will not lead to a
misunderstanding.

\medskip

{\bf Theorem 2}. {\it Let $G$ be an arbitrary quantum space, and
$\cR:E\times F\to G$ an arbitrary strongly completely bounded bioperator.
Then there exists a unique completely bounded operator
$R:E\otimes_hF\to G$ such that the diagram

\medskip

$$
\xymatrix{
E\times F \ar[d]^\vartheta \ar [dr]^\cR\\
E\otimes_hF \ar [r]^R & G      }
$$

\medskip

\noindent is commutative. Moreover, we have $\|R\|_{cb}=\|\cR\|_{scb}$.}

\smallskip

$\tl$ Pure algebra provides a unique linear operator $R$, making our diagram commutative.
The rest follows from Proposition 14. $\tr$

\medskip

Note that, for $u\in\cK E,v\in\cK F$, $\|\vt_s(u,v)\|_h=\|u\bigodot v\|_h\le
\|u\|\|v\|$, and this means that $\vt:E\times F\to E\otimes_hF$ is strongly
completely contractive. Therefore the previous theorem implies

\medskip

{\bf Corollary 3}\quad ("the existence theorem"). {\it The pair
$(E\otimes_hF,\vartheta)$ is the Haagerup
tensor product of quantum spaces $E$ and $F$.}

\medskip

We proceed to the realization of Haagerup tensor product as a projective
tensor product of normed modules.

\medskip

Take, for a moment, an arbitrary normed
algebra $A$, a right normed $A$-module $X$, and a left normed $A$-module $Y$.
Recall that, by definition, the projective tensor product of these modules is
their algebraic module tensor product $X\mma Y$, equipped by the special
semi-norm $\|\cd\|_{mp}$. The latter is the
quotient semi-norm of the projective norm $\|\cd\|_p$ in
$X\otimes Y$ with respect to the canonical quotient map
$\tau:X\otimes Y\to X\mma Y$. (The operator $\tau$ is well defined by
taking $x\otimes y$ to $x\mma y$).
In other words, for $U\in X\mma Y$, $\|U\|_{mp}=
\inf\{\sum_{k=1}^n\|u_k\|\|v_k\|\}$, where the infimum is taken for all
possible representations $U=\sum_{k=1}^nu_k\mma v_k; u_k\in X, v_k\in Y$.
If $X$ and $Y$ are not only one-sided modules but $A$-bimodules,then
$X\mma Y$ also becomes an $A$-bimodule with outer multiplications
well defined by $a\cd(u\otimes_A v):=(a\cd u)\otimes_A v$ and
$(u\mma v)\cd b:=u\mma(v\cd b)$. Moreover, if we consider
$X\otimes Y$ as a tensor product of left and right $A$-modules, $\tau$
becomes a morphism of $A$-bimodules.

In our special context $A=\bb$, $X=\cK E$ and $Y=\cK F$.
As it was mentioned, the strong amplification
$\vt_s:\cK E\times\cK F\to\cK(E\otimes F)$ is a balanced bioperator.
Therefore, by the universal property of the module tensor product, it gives
rise to the linear operator $\odot_\bb:\cK E\mmb\cK F\to\cK(E\otimes F)$,
well defined by $u\mmb v\mapsto\vt_s(u,v)$.

\medskip

{\bf Theorem 3}. {\it The operator $\odot_\bb$ is an isometric $\bb$-bimodule
isomorphism between}
$(\cK E\mmb\cK F,\|\cd\|_{mp})$ and $\cK(E\otimes_hF)$.

\smallskip

$\tl$ Consider the diagram

\medskip

$$
\xymatrix{
(\cK E\otimes\cK F,\|\cd\|_p) \ar[d]^\tau \ar [dr]^\odot\\
(\cK E\mmb\cK F,\|\cd\|_{mp}) \ar [r]^{\odot_\bb} & \cK(E\otimes_h F)      }
$$

\medskip

Recall that $\odot$ and $\tau$ are morphisms of $\bb$-bimodules, and
both of them, by definition of $\|\cd\|_h$ and $\|\cd\|_{mp}$, are
coisometric operators. Since the diagram is obviously commutative,
the operator $\odot_\bb$ also has both properties.
Therefore all what we need is to show that the latter operator is injective.

At first consider the simplest particular case where $E=F=\co$.
Then $\odot_\bb$ is, of course, just the so-called product map
$\pi:\cK\mmb\cK\to\cK:a\mmb b\mapsto ab$. Take
$U\in\cK\mmb\cK$; let it have the
form $U=\sum_{k=1}^na_k\mmb b_k$ for some compact operators $a_k,b_k$.
As is well known (and easy to check), there exist $c,d_k\in\cK$ such
that $a_k=cd_k$ for all $k$ and $c$ is a not a divisor of zero.
Consequently $U=\sum_{k=1}^nc\mmb d_kb_k=
c\mmb(\sum_{k=1}^nd_kb_k)$ and $\pi(U)=c(\sum_{k=1}^nd_kb_k)$.
Therefore, if $\pi(U)=0$, then $\sum_{k=1}^nd_kb_k=0$ and hence
$U=c\mmb0=0$. We see that $\pi$ is injective (and thus it is an
isometric isomorphism).

Now return to general $E$ and $F$. "Changing the order of the relevant
tensor factors", we easily see that the space
$\cK E\mmb\cK F$ coincides, up to a linear isomorphism, with
$(\cK\mmb\cK)\otimes(E\otimes F)$. (To be precise, this isomorphism and its
inverse are defined, by an obvious way, with the help of the 4-linear
operators, acting as $(a,x,b,y)\mapsto(a\mmb b)\otimes(x\otimes y)$
and $(a,b,x,y)\mapsto ax\mmb by; a,b\in\cK,x\in E,y\in F$). Moreover,
under such an identification, the operator $\odot_\bb$ transforms to
$\pi\otimes{\id}:(\cK\mmb\cK)\otimes(E\otimes F)\to\cK(E\otimes F)$.
But we know that $\pi$ is injective. Hence the same is true for
$\pi\otimes{\id}$. The rest is clear. $\tr$

\medskip

If $E$ and $F$ are operator spaces, we can identify $E\otimes F$ with
$E\mHtp F$ and compare Haagerup quantum norm with the standard quantum
norm. The latter will be denoted by $\|\cd\|_{sp}$ and the respective
quantum space by $E\otimes_{sp}F$.

\medskip

{\bf Proposition 18}. {\it Let $E$ and $F$ be operator spaces. Then
we have $\|\cd\|_{sp}\le\|\cd\|_h$.}

\smallskip

$\tl$ By Proposition 12, the bioperator $\vt:E\times F\to E\otimes_{sp}F$ is
strongly completely contractive. Therefore the definition of the Haagerup
tensor product gives, with $\vt$ as $\cR$, that
${\id}:E\otimes_h F\to E\otimes_{sp}F$ is contractive. The rest is clear. $\tr$

\section{The four-named tensor product}

We turn to the explicit construction of the second principal quantum tensor
product. Beginning with the same canonical bioperator $\vt$, now we consider
its weak amplification $\vt_w$. Let us write $u\di v$ instead of
$\vt_w(u,v);u\in\cK E,v\in\cK F$. Of course, this extended "diamond operation"
is well-defined by $ax\di by=(a\di b)(x\otimes y)$ and hence, by
bilinearity, satisfy the identity
$$
(a\di b)\cd(u\di v)\cd(c\di d)=(a\cd u\cd c)\di(b\cd v\cd d).\eqno{(8)}
$$
Let $\di:\cK E\otimes\cK F\to\cK(E\otimes F)$ be the linear operator,
associated with $\vt_w$; it is well-defined by $\di(u\otimes v)=u\di v$.
This operator, contrary to $\odot$, is not bound to be surjective. However,
another, slightly more complicated operator has this attractive property.
We come to this operator after the following observation.

\medskip

{\bf Proposition 19}. {\it Every $a\in\cK$ has the form $b(c\di d)b'$ for some
$b,b',c,d\in\cK$, and even $b(c\di c)b'$ for some $b,b',c\in\cK$.}

\smallskip

$\tl$ Fix an arbitrary orthonormal basis $e_n; n=1,2,...$ in $L$ and put
$e_{m,n}:=i^*(e_m\otimes e_n)$.
It immediately follows from the classical Schmidt Theorem (see,
e.g.,~\cite[Theorem 2.4.1]{baf})
that $a$ has a factorization $IhJ$ where $I,J$ are unitary operators on $L$,
and $h$ is a compact positive operator with eigenvectors $e_{m,n}$. Let
$\lambda_{m,n}$ be the respective eigenvalues.

There exists, of course, a decreasing sequence
$t_n\ge\max\{\lm_{m,n},\lm_{n,m}: m=1,...,n\}$, converging to 0. Then
the double sequence $r_{m,n}:=\sqrt{t_mt_n}$ is not less than $\lm_{m,n}$
and also converges to 0. Therefore $\lm_{m,n}=r_{m,n}s^2_{m,n}$ for some
non-negative $s_{m,n}\le1$. Consider the compact operator
$c'$ well defined by $e_n\mapsto\sqrt{t_n}e_n$ and the bounded operator
$f$ well defined by $e_{m,n}\mapsto s_{m,n}e_{m,n}$. It is easy to check
that $(c'\di c')e_{m,n}=r_{m,n}e_{m,n}$ and hence
$f(c'\di c')f(e_{m,n})=h(e_{m,n})$. Hence $a=If(c'\di c')fJ$. But $c'$
factorizes as $gcg'$ for some $c,g,g'\in\cK$. Consequently $c'\di c'=
(g\di g)(c\di c)(g'\di g')$, and it remains to
put $b:=If(g\di g)$ and $b':=(g'\di g')fJ$. $\tr$

\medskip

Now we introduce the operator ${\uplus}:\cK\otimes\cK E\otimes\cK F\otimes\cK\to
\cK(E\otimes F)$, associated with the 4-linear operator
$(b,u,v,d)\mapsto b\cd(u\di v)\cd d$.

\medskip

{\bf Proposition 20}. {\it The operator ${\uplus}$ is surjective.}

\smallskip

$\tl$ It follows from the previous proposition that an element in
$\cK(E\otimes F)$ of the form $a(x\otimes y)$ is equal to
$b\cd(c\di c)(x\otimes y)\cd d$, that is belongs to the image of
${\uplus}$. It remains to recall that
an arbitrary element in $\cK(E\otimes F)$ is a sum of several elements of the
indicated form. $\tr$

\medskip

Thus $\cK(E\otimes F)$ can be identified with a quotient space of
$\cK\otimes\cK E\otimes\cK F\otimes\cK$. Introduce on the latter space the
projective norm $\|\cd\|_p$ (that is, the projective tensor product of the
four relevant norms), and consider the respective quotient semi-norm on
$\cK(E\otimes F)$. The latter will be denoted by $\|\cd\|_4$. In other words,
our semi-norm is defined by
$$
\|U\|_4:=\inf\{\sum_{k=1}^n\|a_k\|\|u_k\|\|v_k\|\|b_k\|\},\eqno{(9)}
$$
where the infimum is taken for all possible representations of $U$ as
$\sum_{k=1}^na_k\cd(u_k\di v_k)\cd b_k; a_k,b_k\in\cK, u_k\in\cK E, v_k\in\cK F$.
Note that, with respect to $\|\cd\|_p$ and $\|\cd\|_4$, the operator
$\uplus$ is coisometric.

\medskip

{\bf Proposition 21}. {\it The semi-norm $\|\cd\|_4$ in the $\bb$-bimodule
$\cK(E\otimes F)$ satisfies the first axiom of Ruan.}

\smallskip

$\tl$ The proof repeats, with obvious modifications, that of Proposition 13,
and we omit it here. $\tr$

\medskip

{\bf Proposition 22}. {\it Let $G$ be a quantum space, $\cR:E\times F\to G$ a
weakly completely bounded bioperator, $R:E\otimes F\to G$ the associated
linear operator. Then the amplification $R_\ii:\cK(E\otimes F)\to\cK G$ is a
bounded operator with respect to the semi-norm $\|\cd\|_4$ and
the quantum norm on $G$. Moreover, $\|R_\ii\|=\|\cR\|_{wcb}$.}

\smallskip

$\tl$ Consider the 4-linear operator
${\mathcal S}:\cK\times\cK E\times\cK F\times\cK\to\cK G:
(a,u,v,b)\mapsto a\cd\cR_w(u,v)\cd b$.
Since $G$ satisfies $(RI)$, we easily see that the weak
complete boundedness of $\cR$ implies the (usual) boundedness of ${\mathcal S}$,
and $\|{\mathcal S}\|\le\|\cR_w\|$. At the same time we obviously have
$\cR_w(u,v)=\lim_{n\to\ii}{\mathcal S}(P_N,u,v,P_N)\in\cK G$, where $P_N$ is an
approximate identity in $\cK$ consisting of projections. Therefore the
boundedness of ${\mathcal S}$ implies the complete boundedness of $\cR$, and
$\|\cR_w\|\le\|{\mathcal S}\|$. Thus both kinds of the boundedness
are equivalent, and $\|{\mathcal S}\|=\|\cR\|_{wcb}$.

Now consider the diagram

\medskip

$$
\xymatrix{
\cK\otimes\cK E\otimes\cK F\otimes\cK \ar[d]^\uplus \ar [dr]^S\\
\cK(E\otimes F) \ar [r]^{R_\ii} & \cK G      }
$$

\medskip

\noindent where $S$ is the operator, associated with the 4-linear operator
${\mathcal S}$. By the known property of the projective norm, we have
$\|S\|=\|{\mathcal S}\|$. Further, routine calculations with
elementary tensors in $\cK\otimes\cK E\otimes\cK F\otimes\cK$ show that this
diagram is commutative. It follows, taking into account that $\uplus$ is a
coisometric operator, that $\|R_\ii\|=\|S\|$. The rest is clear. $\tr$

\medskip

{\bf Proposition 23}. {\it The estimate $\|\cd\|_h\le\|\cd\|_4$ is valid;
as a corollary, $\|\cd\|_4$ is a norm}.

\smallskip

$\tl$ We know that the canonical bioperator
$\vartheta:E\times F\to E\otimes F$ is strongly completely contractive
with respect to quantum norms on $E,F$ and the (quantum) Haagerup norm on
$E\otimes F$. Hence, by Theorem 1, it is weakly completely contractive with
respect to the same quantum norms. Put it as $\cR$ in the previous
proposition. In this situation $R$ is, of course, the identity operator in
$E\otimes F$, and $R_\ii$  is the identity operator from
$(\cK(E\otimes F),\|\cd\|_4)$ onto $(\cK(E\otimes F),\|\cd\|_h)$. By the same
proposition, $\|R_\ii\|\le1$. The rest is clear. $\tr$

\medskip

>From now on, we shall call $\|\cdot\|_4$ the {\it four-named
norm}.\footnote{"Blecher/Paulsen-Effros/Ruan norm"; cf. above.}

The role of the following observation concerning the introduced norm is
similar to that of Proposition 16 for the Haagerup norm.

\medskip

{\bf Proposition 24}. {\it Every $U\in\cK(E\otimes F)$ can be represented as}
(a "single rigged diamond")
$$
a\cd(u\di v)\cd b
$$
where $a,b\in\cK, u\in\cK E, v\in\cK F$. {\it In more detail, if
$U=\sum_{k=1}^na_k\cd(u_k\di v_k)\cd b_k;
a_k,b_k\in\cK, u_k\in\cK E, v_k\in\cK F$, and $S_1,...,S_n$ are some isometric
operators with pairwise orthogonal images, then, to obtain such a
representation, one can take
$a:=\sum_{k=1}^na_k(S_k^*\di S_k^*), u:=\sum_{k=1}^nS_k\cd u_k\cd S_k^*,
v:=\sum_{k=1}^nS_k\cd v_k\cd S_k^*$ , and
$b:=\sum_{k=1}^n(S_k\di S_k)b_k$. Finally, we have
$$
\|U\|_4:=\inf\{\|a\|\|u\|\|v\|\|b\|\},\eqno{(10)}
$$
where the infimum is taken for all possible representations of $U$ in the
indicated form.}

\smallskip

$\tl$ Recall that $S_k^*S_l=\de^k_l$.
Therefore the routine calculation using the equalities (3) and
(8) shows that $U$ indeed has the desired representation.

To obtain the desired equality for $\|U\|_4$, take $\e>0$. By (9), there
exists a representation
$U=\sum_{k=1}^na_k\cd(u_k\Diamond v_k)\cd b_k$ such that
$\sum_{k=1}^n\|a_k\|\|u_k\|\|v_k\|\|b_k\|<\|U\|_4+\e$. Using, if necessary,
scalar multiples,
we have a right to assume that $\|a_k\|=\|b_k\|$ and $\|u_k\|=\|v_k\|=1$
for all $k$, and thus $\sum_{k=1}^n\|a_k\|^2<\|U\|_4+\e$. Now take the
representation of $U$ as a "single rigged diamond", indicated above.
We see that $u$ is the sum of several elements of norm 1 with the pairwise
orthogonal supports, namely $S_kS_k^*$, and the same is true for $v$.
Therefore, by $(RII)$, we have $\|u\|=\|v\|=1$. Finally, the operator
$C^*$-identity, together with the formula (3), gives
$$
\|a\|=\left|\left|\left[\sum_{k=1}^na_k(S_k^*\di S_k^*)\right]
\left[\sum_{l=1}^n(S_l\di S_l)a_l^*\right]\right|\right|^{\frac{1}{2}}=
$$
$$
\left|\left|\sum_{k,l=1}^na_k[(S_k^*S_l)\di(S_k^*S_l)]
a_l^*]\right|\right|^{\frac{1}{2}}=
\left|\left|\sum_{k=1}^na_ka_k^*\right|\right|^{\frac{1}{2}}
\le\left(\sum_{k=1}^n\|a_k\|^2\right)^{\frac{1}{2}}.
$$

Similar calculations, combined with $\|a_k\|=\|b_k\|$, give the same
estimation for $\|b\|$. Consequently, $\|a\|\|u\|\|v\|\|b\|\le
\sum_{k=1}^n\|a_k\|^2<\|U\|_4+\e$. Since $\e>0$ is arbitrary, this implies
that the infimum in (10) is not bigger than
$\|U\|_4$. The inverse inequality is obvious. $\tr$

\medskip

{\bf Proposition 25}. {\it The four-named norm on $\cK(E\otimes F)$ satisfies
the second axiom of Ruan.}

\smallskip

$\tl$ Let $U,V\in\cK(E\otimes F)$ have orthogonal supports $P_1$ and $P_2$.
Again (as in the proof of Proposition 17) we have a right to assume that
$\|U\|_4>\|V\|_4$ and that
the images of $P_1$ and $P_2$ are infinite-dimensional.

Using the previous proposition, take an arbitrary $\e$ with
$0<\e<\|U\|_4-\|V\|_4$ and representations
$U=a_1\cd(u_1\Diamond v_1)\cd b_1$ and $V=a_2\cd(u_2\Diamond v_2)\cd b_2$
such that $\|a_1\|\|u_1\|\|v_1\|\|b_1\|<\|U\|_4+\e$ and
$\|a_2\|\|u_2\|\|v_2\|\|b_2\|<\|V\|_4+\e$. Of course, we
can assume that $a_k=P_ka_k, b_k=b_kP_k, \|u_k\|=\|v_k\|=1; k=1,2$, and also
$\|a_1\|\ge\|a_2\|$ and $\|b_1\|\ge\|b_2\|$. In particular, we have
$\|a_1\|\|b_1\|<\|U\|_4+\e$.

Now take isometric operators $S_k;k=1,2$ with $S_kS^*_k=P_k$ and put
$a=a_1(S_1^*\di S_1^*)+a_2(S_2^*\di S_2^*),
u=S_1\cd u_1\cd S_1^*+S_2\cd u_2\cd S_2^*,
v=S_1\cd v_1\cd S_1^*+S_2\cd v_2\cd S_2^*$,
and $b=(S_1\di S_1)b_1+(S_2\di S_2)b_2$.
By the part of Proposition 24, presenting the r\'ecipe of a "single
diamond", we have
$$
U+V=a\cd(u\Diamond v)\cd b.
$$
It follows that $\|U+V\|_4\le\|a\|\|u\|\|v\|\|b\|$. Since the elements
$S_k\cd u_k\cd S_k^*;k=1,2$ have supports $P_k$, the axioms $(RII)$ and
then $(RI)$ give $\|u\|=\max\{\|S_k\cd u_k\cd S_k^*\|;k=1,2\}=1$.
Similarly $\|v\|=1$. Finally, the operators $a_k(S_k^*\di S_k^*);k=1,2$
have orthogonal left supports $P_k;k=1,2$ and, by (3), orthogonal right
supports $P_k\di P_k; k=1,2$.

As is well known in operator theory, these properties of the summands of
$a$ imply that $\|a\|=\max\{\|a_k(S_k^*\di S_k)\|;k=1,2\}$. Hence, by
Proposition 1,
we have $\|a\|=\max\{\|a_k\|;k=1,2\}=\|a_1\|$. Similarly $\|b\|=\|b_1\|$.
Therefore $\|U+V\|_4\le\|a_1\|\|b_1\|<\|U\|_4+\e$. Since $\e>0$ is arbitrary,
we have $\|U+V\|_4\le\|U\|_4$. Because of $P_1\cd(U+V)=U$, the inverse
inequality follows from $(RI)$. The rest is clear. $\tr$

\smallskip

Combining Propositions 21 and 25, we obtain

\medskip

{\bf Corollary 4}. {\it The four-named norm is a quantum norm on
$E\otimes F$.}

\medskip

We denote the constructed quantum space by $E\otimes_4F$. The same symbol
("the sign of the four") will be denote the underlying normed space;
this will not lead to a misunderstanding.

\medskip

{\bf Theorem 4}. {\it Let $G$ be an arbitrary quantum space, and
$\cR:E\times F\to G$ an arbitrary weakly completely bounded bioperator. Then
there exists a unique completely bounded operator
$R:E\otimes_4F\to G$ such that the diagram

\medskip

$$
\xymatrix{
E\times F \ar[d]^\vartheta \ar [dr]^\cR\\
E\otimes_4F \ar [r]^R & G      }
$$

\medskip

\noindent is commutative. Moreover, we have $\|R\|_{cb}=\|\cR\|_{wcb}$.}

\smallskip

$\tl$ The argument of Theorem 2 works, with Proposition 22 replacing
Proposition 14. $\tr$

\medskip

Note that, for $u\in\cK E,v\in\cK F$, we have $\|\vt_w(u,v)\|_4=
\|u\Diamond v\|_4\le\|u\|\|v\|$. This means that
$\vt:E\times F\to E\otimes_4F$ is weakly completely contractive. Therefore
the previous theorem implies

\medskip

{\bf Corollary 5}\quad ("the existence theorem"). {\it The pair
$(E\otimes_4F,\vartheta)$ is the
four-named tensor product of quantum spaces $E$ and $F$.}

\section{Examples}

Again, to get instructive illustrations, we turn to our well-beloved column
and raw Hilbertians. Recall the identifications of $H_c$ with $\bb(\co,H)$
and of $H_r$ with $\bb({\overline H},\co)$. In what follows, the symbols
$H_c$ and $H_r$ denote, depending on the context, the respective standard
quantum spaces or their underlying normed spaces; this will not lead to a
confusion.

\medskip

{\bf Proposition 26}. {\it Let $H$ be a Hilbert space, $E$ an arbitrary
operator
space. Then, up to complete isometric isomorphisms, $H_c\otimes_h E=
H_c\otimes_{sp} E$ and $E\otimes_h H_r=E\otimes_{sp}H_r$. More precisely,
the identity operators ${\id}:H_c\otimes_h E\to H_c\otimes_{sp}E$ and
${\id}:E\otimes_h H_r\to E\otimes_{sp}H_r$ are complete isometric
isomorphisms.}

\smallskip

$\tl$ We already know, by Proposition 18, that both identity operators are
completely contractive. Therefore our task is to show that their
amplifications do not decrease norms.

Consider the "column" case. Take $U\in\cK(H\otimes E)$. Identifying the
latter space with $H\otimes(\cK E)$ and using the first part of
Proposition 5, we can represent $U$ as $\sum_{k=1}^ne_k\otimes u_k$ with
$e_k$ as in that proposition and $u_k\in\cK E$. (Of course, it does not matter
that we have now the order of tensor factors different from that in the
cited proposition).

Fix, for a time, an element of the form $\omega\in\cK H_c$ indicated in
(2) together with the relevant projection $P$ and partial isometries
$q_k; k=1,...,n$. Put $u:=\sum_{k=1}^nq_k\cd u_k$.
Using that $q_k^*q_l=\de^k_lP$, we have
$$
\omega\odot u=\sum_{k,l=1}^nq_k^*e_k\odot q_l\cd u_l=
\sum_{k,l=1}^nq_k^*q_le_k\odot u_l=\sum_{k=1}^nPe_k\odot u_k.
$$
Assuming that $u_k$ is an elementary tensor in $\cK E$, one can easily check
that $Pe_k\odot u_k\in\cK(H\otimes E)$ is exactly $P\cd(e_k\otimes u_k)$.
Then, by bilinearity, the same is true in the general case. It obviously
follows that $\omega\odot u=P\cd U$.

By Proposition 6, we have $\|\omega\|=1$. Further, elements $q_k\cd u_k$
live in the operator space $\cK\mHtp E$. Therefore, by $C^*$-identity we have
$$
\|u\|_{sp}=\left|\left|\left(\sum_{k=1}^nq_k\cd u_k\right)^*
\left(\sum_{l=1}^nq_l\cd u_l\right)\right|\right|^\frac{1}{2}=
\left|\left|\sum_{k,l=1}^n[(q_k\mHtp{\id})u_k]^*
[(q_l\mHtp{\id})u_l]\right|\right|^\frac{1}{2}=
$$
$$
\left|\left|\sum_{k,l=1}^nu_k^*[q_k^*q_l\mHtp
{\id}]u_l\right|\right|^\frac{1}{2}=
\left|\left|\sum_{k=1}^nu_k^*(P\mHtp
{\id})u_k\right|\right|^\frac{1}{2}.
$$
Consequently, by the definition of the Haagerup norm, we have
$$
\|P\cd U\|_h\le\|\omega\|\|u\|\le
\left|\left|\sum_{k=1}^nu_k^*(P\mHtp
{\id})u_k\right|\right|^\frac{1}{2}.
$$

Now consider a sequence $P_N$ of finite-dimensional projections, serving as
an approximate identity in $\cK$. Taking elementary tensors in $\cK\mHtp E$
and using the bilinearity, we see that $\sum_{k=1}^nu_k^*(P_N\mHtp{\id})u_k$
converges, with respect to the operator norm, to $\sum_{k=1}^nu_k^*u_k$.
At the same time, of course, $P_N\cd U$ converges to
$U$ in $\cK(H_c\otimes_hE)$. Combined with the obtained inequality,
both things give the estimate
$$
\|U\|_h\le\|\omega\|\|u\|\le\|\sum_{k=1}^n u_k^*u_k\|^\frac{1}{2}.
$$
Now turn to the norm of $U$ as of an element of $\cK\mHtp(H_c\mHtp E)$ or,
equivalently, of $H_c\mHtp(\cK\mHtp E)$. This time, by virtue of Proposition 5,
we have the exact equality $\|U\|_{sp}=
\|\sum_{k=1}^n u_k^*u_k\|^\frac{1}{2}$.

This ends the proof in the "column" case. The similar argument, with the
obvious modifications, works in the "raw" case. $\tr$

\medskip

Here is an illuminating particular case. Denote by ${\mathcal F}(H)$ the space
of bounded finite-dimensional operators on $H$, endowed by the operator norm
and the standard quantization.

\medskip

{\bf Proposition 27}. Let $H$ be a Hilbert space. Then, up to a
complete isometric isomorphism, $H_c\otimes_h{\overline H}_r={\mathcal F}(H)$.

\smallskip

$\tl$ By virtue of the previous proposition, it is sufficient to establish a
complete isometric isomorphism between standard quantum spaces
$H_c\otimes_{sp}{\overline H}_r$ and ${\mathcal F}(H)$.

The space $H_c\otimes_{sp}{\overline H}_r$ or, otherwise,
$\bb(\co,H)\mHtp\bb(H,\co)$, is a subspace in $\bb(\co\mHtp H,H\mHtp\co)$.
The latter, because of the identification of $\co\mHtp H$ and $H\mHtp\co$
with $H$, coincides with $\bb(H)$. It is easy to see that the resulting
isometric embedding of $H_c\otimes_{sp}{\overline H}_r$ into $\bb(H)$ takes
an elementary tensor $x\mHtp y$ to the rank 1
operator $x\bigcirc y$. Obviously, the image of this embedding is
${\mathcal F}(H)$.
Denote by $I:H_c\otimes_{sp}{\overline H}_r\to{\mathcal F}(H)$ the respective
corestriction. Now it is sufficient for us to show
that its amplification $I_\ii$ is an isometric isomorphism.

Since our quantum spaces are standard, $I_\ii$ acts between
the subspace $\cK\mHtp H_c\mHtp{\overline H}_r$ in
$\bb(L\mHtp\co\mHtp H,L\mHtp H\mHtp\co)$ and the subspace
$\cK\mHtp{\mathcal F}(H)$ in $\bb(L\mHtp H)$,
and it is uniquely determined by taking
$a\mHtp(x\mHtp y);a\in\cK$ to
$a\mHtp(x\bigcirc y)$. From this we easily see that $I_\ii$ is a
birestriction of a certain isometric isomorphism between  these bigger
spaces. The latter is generated by the
natural identification of $L\mHtp\co\mHtp H$ and $L\mHtp H\mHtp\co$ with
$L\mHtp H$. Therefore
$I_\ii$ is itself an isometric isomorphism. The rest is clear. $\tr$

\medskip

We have described what happens if the left factor in the Haagerup tensor
product is a column
Hilbertian. But what if we put this Hilbertian on the right?

\medskip

{\bf Proposition 28}. {\it Let $H$ be a Hilbert space, $E$ an arbitrary
quantum  space. Then, up to complete isometric isomorphisms,
$E\otimes_h H_c=E\otimes_4H_c$ and $H_r\otimes_h E=H_r\otimes_4E$. More
precisely, the identity operators ${\id}:E\otimes_4H_c\to E\otimes_h H_c$
and ${\id}:H_r\otimes_4E\to H_r\otimes_hE$ are
complete isometric isomorphisms.}

\smallskip

$\tl$ By virtue of Proposition 23, our task is only to show that the
amplifications of our identity operators do not decrease norms.

Consider the "column" case. Take $U\in\cK(E\otimes H)$ and, using
Proposition 16, represent it as a single Effros symbol
$u\odot v;u\in\cK E,v\in\cK H$. Further, using Proposition 5, represent
$v$ as $\sum_{k=1}^na_ke_k;a_k\in\cK$, where $e_k$ is an orthonormal system
in $L$.

Taking some $P$ and $q_k$, consider an element of the form $\omega$ as
indicated in (2). Put
$b:=\sum_{k=1}^na_k\di q_k\in\cK$. Then, using the equality (8), we have
$$
(u\di\omega)\cd b=\sum_{k,l=1}^n(u\di q_k^*e_k)\cd(a_l\di q_l)=
\sum_{k,l=1}^n(u\cd a_l)\di(q_k^*q_l)e_k=\sum_{k=1}^n(u\cd a_k)\di Pe_k.
$$
If $u$ is an elementary tensor, one can easily verify that
$(u\cd a_k)\di Pe_k=\\(u\odot a_ke_k)\di P$. Hence, by bilinearity, the same
is true for the general $u\in\cK E$. Therefore we have
$$
(u\di\omega)\cd b=\sum_{k=1}^n(u\cd a_k)\di Pe_k=
\sum_{k=1}^n(u\odot a_ke_k)\di P=(u\odot v)\di P=U\di P.
$$
>From this, combining the expression of the four-named norm by (10)
with Propositions 8 and 6,
we see that $\|U\|_4=\|U\di P\|_4\le\|u\|\|\omega\|\|b\|=\|u\|\|b\|$.
But the $C^*$-identity, together with Propositions 8 and 5, gives
$\|b\|=\|v\|$. Therefore,
taking all possible representations of $U$ as single Effros symbols and
using (7), we complete the proof in the "column" case. The similar argument,
with the obvious modifications, works in the "raw" case. $\tr$

\medskip

The following important observation illustrates both Propositions 26 and 28.
Let $H$ be a Hilbert space. Denote by $(H\otimes H)_c$ (respectively,
$(H\otimes H)_r$ the algebraic tensor square of $H$, considered as a quantum
subspace of the column Hilbertian $(H\mHtp H)_c$ (respectively, raw
Hilbertian $(H\mHtp H)_r$).

\medskip

{\bf Proposition 29}.
{\it Up to complete isometric isomorphisms, $H_c\otimes_4 H_c=\\H_c\otimes_hH_c=
H_c\otimes_{sp}H_c=(H\otimes H)_c$ and $H_r\otimes_4H_r=H_r\otimes_hH_r=
H_r\otimes_{sp}H_r=\\(H\otimes H)_r$.}

\smallskip

$\tl$ Because of the analogy between the "column" and the "raw" cases, it is sufficient to restrict
ourselves with the first chain of equalities. In this chain, the first two equalities follow from
the mentioned propositions. We proceed to the third equality.

Recall that $H_c=\bb(\co,H)$ and consider the linear isomorphism
$I:H_c\mHtp H_c\to(H\otimes H)_c$,
coinciding, after the respective identifications, with the identity
operator on $H\otimes H$. Since we deal with standard quantum
spaces, its amplification $I_\ii$ acts between $\cK\mHtp H_c\mHtp H_c$
and $\cK\mHtp(H\otimes H)_c$. This is obviously a birestriction of a certain
isometric isomorphism. The latter, if we want to be meticulous, acts
between $\bb(L\mHtp\co\mHtp\co,L\mHtp H\mHtp H)$ and
$\bb(L\mHtp\co,L\mHtp H\mHtp H)$ and is
generated by the natural identification of $L\mHtp(\co\mHtp\co)$ with
$L\mHtp\co$. Therefore
$I_\ii$ is itself an isometric isomorphism. The rest is clear. $\tr$

\bibliographystyle{amsalpha}

\ed